%rna.tex: a Plain TeX file by AJ Bu, Manuel Kauers and  and Doron Zeilberger
%Statistical Analysis of Hairpins and BasePairs in RNA Secondary Structures

%begin macros

\baselineskip=14pt
\parskip=10pt

\magnification=\magstephalf

\def\P{{\cal P}}

\def\1{{\overline{1}}}
\def\2{{\overline{2}}}
\def\E{{\rm E}}
\def\Var{{\rm Var}}
\def\Corr{{\rm Corr}}
\def\Cov{{\rm Cov}}
\parindent=0pt
\overfullrule=0in

\def\frac#1#2{{#1 \over #2}}
%\headline={\rm  \ifodd\pageno  \RightHead  \else  \LeftHead  \fi}
%\def\RightHead{\centerline{
%Title
%}}
%\def\LeftHead{ \centerline{Doron Zeilberger}}
%end macros
\centerline
{\bf
Statistical Analysis of Hairpins and BasePairs in RNA Secondary Structures
}
\bigskip
\centerline
{\it AJ BU, Manuel KAUERS, and Doron ZEILBERGER}
\bigskip

\qquad {\it In memory of Kequan Ding (1955-2025)}

{\bf Abstract}: We derive precise asymptotic expressions for the
expectations, variances, covariance, and quite a few further mixed moments for 
the number of hairpins and the number of basepairs in RNA secondary structures, and give convincing
evidence that the central-scaled distribution of the pair of random variables
(hairpins, basepairs) tends in distribution to the bi-variate normal distribution
with correlation  $\sqrt{5 \sqrt{5} -11}/2= 0.2123322205\dots$

{\bf Dedication:} This article is in fond memory of Kequan Ding, a very deep, original, and versatile mathematician, whose
many interests included RNA secondary structures [GD].

{\bf Introduction}: In  a seminal paper [W], Michael Waterman presented a combinatorial framework for RNA secondary structures. 
A secondary structure is a vertex-labeled graph on $n$ vertices obeying certain conditions. See [W], [DG], and [CRU].
These have several {\it structure elements}, the most important one being {\it basepair} and another one
is {\it hairpin}. We refer to [CRU] (readily available on-line) for the definitions.

The following theorem is due to Waterman [W]:

{\bf Theorem 1}: Let $a(n)$ be the number of RNA secondary structures on $n$ vertices. Then
$$
\sum_{n=0}^{\infty} a(n) X^n=
\frac{-2 X^{4}+X^{3}+2 X -1+\sqrt{X^{6}-4 X^{5}+4 X^{4}-2 X^{3}+4 X^{2}-4 X +1}}{2 X^{2} \left(X -1\right)} \quad .
$$

Using symbolic computation, one of us [B], proved the following generalization.

{\bf Theorem 2}: Let $A(n,i,j)$ be the number of RNA secondary structures on $n$ vertices, with exactly
$i$ hairpins and $j$ basepairs. Let's define the tri-variate generating function $F(X;x,z)$
$$
F(X;x,z):=\sum_{n=0}^{\infty} \sum_{i=0}^{n} \sum_{j=0}^{n} A(n,i,j) \, X^n \, x^i \, z^j \,  \quad,
$$
viewed as a formal power series in $X$ with coefficients that are polynomials in $x$ and $z$. We have the
following {\it hairy} formula
$$
F(X;x,z) \, = \,
\frac{-2 X^{4} z +X^{3} x z +X^{2} z -X^{2}+2 X -1+\sqrt{R}}{2 X^{2} z \left(X -1\right)} \quad ,
$$
where
$$
\eqalign{
R&=X^{6} x^{2} z^{2}-2 X^{5} x \,z^{2}-2 X^{5} x z +4 X^{4} x z +X^{4} z^{2}-2 X^{4} z -2 X^{3} x z\cr
 &\qquad+X^{4}+4 X^{3} z -4 X^{3}-2 X^{2} z +6 X^{2}-4 X +1 \quad .
}
$$

In this {\it methodological} article (the method is much more important than the actual results!), we will establish the following theorem.

{\bf Theorem 3}: For a secondary structure $s$ on $n$ vertices,
let $X_n(s)$ and $Z_n(s)$ be the number of hairpins and number of basepairs, respectively, of $s$.
We have the following asymptotic expressions
$$
\E[X_n]=(1- \frac{2}{\sqrt{5}})\, n \cdot \left (1 + O(\frac{1}{n})  \right ) \, = \,
  (0.1055728092\dots) \cdot \, n \cdot \left (1 + O(\frac{1}{n})  \right ) \quad .
$$
$$
\E[Z_n]=\frac{5-\sqrt{5}}{10}\, n \cdot \left (1 +  O(\frac{1}{n})   \right ) \, = \,
(0.2763932023\dots) \cdot \, n \cdot \left (1 +  O(\frac{1}{n})   \right )  \quad .
$$
$$
\Var(X_n)=(2- \frac{22}{5 \sqrt{5}})\, n \cdot \left (1 +  O(\frac{1}{n})  \right) \, = \,
(0.032260180\dots)\, \cdot n \cdot \left (1 +  O(\frac{1}{n})  \right) \quad .
$$
$$
\Var(Z_n)= \frac{1}{10 \sqrt{5}}\, n \cdot \left (1 +  O(\frac{1}{n}) \right ) \, = \,
(0.04472135954\dots)\, \cdot n \cdot \left (1 +  O(\frac{1}{n})  \right) \quad .
$$
$$
\Corr(X_n,Z_n)=\frac{\sqrt{5 \sqrt{5} -11}}{2} \cdot \left (1+  O(\frac{1}{n})  \right ) \, = \,
(0.2123322161\dots) \cdot \left (1+  O(\frac{1}{n})  \right ) \quad .
$$

More precisely, we have:

$$
\eqalign{
\E[X_n]&= (1-\frac{2}{5}\sqrt{5})n \cdot \Bigl(1+(\frac{7}{4}+\frac{11}{20}\sqrt{5})n^{-1}+(\frac{243}{160}+\frac{99}{160}\sqrt{5})n^{-2}+(\frac{339}{160}+\frac{1029}{800}\sqrt{5})n^{-3}\cr
      &+(\frac{11917563}{819200}+\frac{5400687}{819200}\sqrt{5})n^{-4}+{\rm O}(n^{-5})\Bigr)\cr
\E[Z_n]&= (\frac{1}{2}-\frac{1}{10}\sqrt{5})n \cdot \Bigl(1+(-\frac{5}{8}-\frac{13}{40}\sqrt{5})n^{-1}+(\frac{21}{320}+\frac{21}{320}\sqrt{5})n^{-2}\cr
      &+(-\frac{93}{320}-\frac{177}{1600}\sqrt{5})n^{-3}+(-\frac{13887249}{1638400}-\frac{6272793}{1638400}\sqrt{5})n^{-4}+{\rm O}(n^{-5})\Bigr)\cr
\Var[X_n] &= (2-\frac{22}{25}\sqrt{5})n \cdot \Bigl(1+(\frac{1}{16}+\frac{23}{80}\sqrt{5})n^{-1}+(-\frac{651}{160}-\frac{177}{64}\sqrt{5})n^{-2}\cr
     &+(-\frac{1208783}{81920}-\frac{3216693}{409600}\sqrt{5})n^{-3}+{\rm O}(n^{-4})\Bigr)\cr
\Var[Z_n] &= (0+\frac{1}{50}\sqrt{5})n \cdot \Bigl(1+(1+\frac{1}{10}\sqrt{5})n^{-1}-\frac{261}{160}n^{-2}\cr
     &+(-\frac{27179}{2560}-\frac{915879}{204800}\sqrt{5})n^{-3}+{\rm O}(n^{-4})\Bigr)\cr
\Corr[X_n,Z_n] &= \frac{1}{2} \sqrt{5\sqrt{5}-11}\cdot \Bigl(1+(\frac{15}{32}+\frac{13}{32}\sqrt{5})n^{-1}+(\frac{4469}{1024}+\frac{1345}{1024}\sqrt{5})n^{-2}\cr
     &+(\frac{1766711}{32768}+\frac{801983}{32768}\sqrt{5})n^{-3}+{\rm O}(n^{-4})\Bigr)
}
$$

We have  lots of evidence that the following conjecture is true, and one of us (DZ) is offering a donation of \$100 to the OEIS for its rigorous proof.

{\bf Conjecture}: The centralized-scaled version of the pair $(X_n,Z_n)$, namely
$$
 \left (  \frac{ X_n-\E[X_n] }{\sqrt{\Var(X_n)}}, \frac{Z_n-\E[Z_n]}{\sqrt{\Var(Z_n)}}  \right ) \quad ,
$$
tends, in distribution, to the bi-variate normal distribution with covariance
$c=\frac{\sqrt{5 \sqrt{5} -11}}{2}$, whose probability density function (pdf) is:
$$
\frac{{ e}^{-\frac{1}{2} x^{2}-\frac{1}{2} y^{2}+c x y} \sqrt{1-c^{2}}}{2 \pi} \quad .
$$

{\bf Reminders about using Symbolic Computations to do Statistics}

Given a finite set $S$, equipped with some random variable, $R$, $s \rightarrow R(s)$, the {\it expectation}, $\mu$, is, of course
$$
\mu=\E[R]:=\frac{1}{|S|}\, \left (\sum_{s \in S} R(s) \right ) \quad,
$$
where, as usual, $|S|$ denotes the number of elements in $S$.
More generally, the $i$-th moment of $R$, is the average of the $i$-th power:
$$
\E[R^i]:=\frac{1}{|S|}\, \left (\sum_{s \in S} (R(s))^i \right ) \quad .
$$

More useful, are the {\bf central moments}
$$
\E[(R-\mu)^i]:=\frac{1}{|S|}\, \left (\sum_{s \in S} (R(s)-\mu)^i \right ) \quad ,
$$
the most important being the {\it second}, aka {\it variance}.
Using the binomial theorem, once you know $\mu$, you can derive the central moments from the straight moments.

A convenient way to record the data regarding the random variable $R$ defined on a finite (and even infinite) set, is via
its {\it weight-enumerator}, aka {\it generating function}, pick any (formal) variable name, say $x$, and define:
$$
f_R(x):=\sum_{s \in S} x^{R(s)} \quad .
$$
In terms of $f_R(x)$ we have
$$
\mu=\E[R]=f_R'(1)= (x \frac{d}{dx}) f_R(x) \Bigl \vert_{x=1} \quad .
$$
More generally, the $i$-th moment is
$$
m_i=\E[R^i]=(x \frac{d}{dx})^i f_R(x) \Bigl \vert_{x=1} \quad .
$$
Once you know the expectation, $\mu$, and the {\it variance}, (aka second central moment), $\sigma^2=E[(R-\mu)^2]$, (where
$\sigma$, the square-root of the variance, is called the {\it standard deviation}), you know the two most important
numbers, what it is `on average', and how spread-out it is. But it is still nice to know more, and then
one forms the {\it scaled-central distribution}
$$
\bar{R} (s) := \frac{R(s)-\mu}{\sigma} \quad,
$$
whose expectation is $0$ and variance is~$1$. Its higher moments are expressible in terms of the moments of the original
random variable~$R$, and hence a knowledge of the weight-enumerator, $f_R(x)$, is useful.

In combinatorics, we often have an {\it infinite} family of sets $S_n$, naturally indexed by non-negative integers, and  a natural
random variable defined on them, with a `nice' sequence of weight-enumerators $f_n(x)$. For example for the set of subsets
of $\{1,\dots, n\}$ with {\it weight} `cardinality', we have the closed-form expression, $f_n(x)=(1+x)^n$.
If not so lucky, we often have, in combinatorics, that the {\bf grand generating function},
$$
F(X;x):=\sum_{n=0}^{\infty} f_n(x) X^n \quad,
$$
has some explicit expression in $X$ and $x$, or otherwise, in many combinatorial families, satisfies some
{\it algebraic equation}
$$
\P (F,X,x) \, = \, 0 \quad,
$$
where $\P$ is some  specific  polynomial in its arguments.

If we have an `explicit' (or even `implicit') expression for the grand-generating function, that is a certain formal power series
in $X$ with coefficients that are polynomials in $x$, we can get generating functions for the numerical sequences of the moments.

Of course, we also need to handle the cardinality of the family of sets, let's call is $a_0(n)$, whose
generating function is
$$
\sum_{n=0}^{\infty} a_0(n) \, X^n =F(X;1) \quad.
$$
Defining
$$
\sum_{n=0}^{\infty} a_i(n) \, X^n = (x \frac{d}{dx})^i F(X;x) \Bigl \vert_{x=1} \quad,
$$
the $i$-th  moment is then $\frac{a_i(n)}{a_0(n)}$.

Alas, it is generally not possible to have closed-form expressions for $a_i(n)$ and not even for $a_0(n)$, but one
can have the next-best thing, {\it asymptotics}. If the grand-generating function is {\it algebraic}, and
especially if it is quadratic, it is relatively easy to find the {\it leading} asymptotics of the
coefficients (See [KP], Section 6.5, and the classic [FS]). Alas, for the variance and higher central moments,
because of cancellations, we need the more challenging problem of having higher-order asymptotics for these quantities.

The good news is that at least for quadratic algebraic formal power series, this can be done effectively.

{\bf Pairs of random variables}

Often our combinatorial set is equipped with several random variables, and in addition to wanting to find the asymptotics for the individual
expectations, means, and higher moments, we'd love to know how they {\it interact}. If we have two random variables, $R_1(s)$, and $R_2(s)$,
then the most important quantity is the {\it covariance}
$$ 
\Cov(R_1,R_2) := \E[(R_1-\mu_1)(R_2-\mu_2)]= \E[R_1 R_2] - \E[R_1]\,\E[R_2] \quad,
$$
and more generally, the {\it central mixed moments}
$$
m_{i,j} := \E[(R_1-\mu_1)^i(R_2-\mu_2)^j] \quad,
$$
that by using the binomial theorem can be expressed in terms of the moments $\{\E[R_1^i R_2^j]\}$.
Even more important than the covariance is the {\it correlation} (that is always between $-1$ and $1$):
$$
\Corr(R_1,R_2):=\frac{\Cov(R_1,R_2)}{\sqrt{\Var(R_1)\,\Var(R_2)}} \quad .
$$

If you define the bi-variate weight-enumerator
$$
f(x_1,x_2):= \sum_{s \in S} x_1^{R_1(s)} x_2^{R_2(s)} \quad,
$$
then again all the relevant quantities can be obtained by repeatedly applying $x_1 \frac{d}{dx_1}$ and/or $x_2 \frac{d}{dx_2}$ and
at the end of the day plugging-in $x_1=1,x_2=1$.

Going back to combinatorial families. You have an infinite sequence of related sets, e.g. RNA secondary structures on $n$ vertices, and two
natural random variables defined on them, e.g. `number of hairpins' and `number of basepairs', and you have the grand-generating function
$$
F(X;x_1,x_2)= \sum_{n=0}^{\infty} f_n(x_1,x_2) \, X^n \quad.
$$
By applying $(x_1 \frac{d}{dx_1})^i (x_2 \frac{d}{dx_2})^j$, and then plugging-in $x_1=1,x_2=1$,
we get formal powers series in $X$, with {\it numerical} coefficients, for the $(i,j)$-mixed moment, and then and use the algorithms to find higher-order asymptotics for the
needed quantities, in particular for the covariance, and from there, the correlation.

{\bf Let's Keep it Simple: Finding the asymptotics approximately }

Once {\it symbolic} computation gave us explicit expressions for the relevant generating functions, we can get lots and lots of coefficients, and since we know from {\it general nonsense},
(we are dealing with the {\it Sch\"utzenberger ansatz} (algebraic formal power series)),
that the asymptotics, in each case, must be of the form
$$
\mu^n n^\theta \cdot (C_0 + \frac{C_1}{n}+ \frac{C_2}{n^2}+ \dots ) \quad,
$$
for {\it some} numbers $\mu$, $\theta$, and $C_0,C_1,C_2, \dots$.

Both $\mu$ and $\theta$ are easily found, and even $C_0$ is relatively easy (see [KP], [FS]), but $C_1$, $C_2$, $\dots$, are harder.
But by using {\it numerics} you can do {\it regression} and get good estimates, that for all practical purposes are {\it good enough}.

(Please note that this $\mu$ has nothing to do with the former $\mu$. Here we bow to two traditions. In statistics, $\mu$ denotes the expectation (alias mean, alias average),
while in Statistical Physics, $\mu$ denotes the so-called {\it connective constant}, which is the exponential rate of growth.
Let us also mention that  $\theta$ is called the {\it critical exponent}).

This is implemented in the Maple package {\tt RNAstat.txt} available from:

{\tt http://sites.math.rutgers.edu/\~{}zeilberg/tokhniot/RNAstat.txt} .

For a detailed numeric statistical analysis see the output file:

{\tt http://sites.math.rutgers.edu/\~{}zeilberg/tokhniot/oRNAstat1.txt} \quad .

But it is still nice to get the relevant constants {\it exactly}, in terms of algebraic numbers. While these can sometimes be {\it guessed} using, e.g.
Maple's {\it identify} command (based on PSLQ and LLL algorithms), these are not rigorous, and even then, it is often not even possible to get a reasonable guess.

In the next section we will describe how we proved Theorem~3.

{\bf Finding the Exact Asymptotics}

Every algebraic power series satisfies a linear differential equation with
polynomial coefficients, and every such differential equation can be easily
translated into a linear recurrence equation with polynomial coefficients for
the coefficient sequence of the series.

Given such a recurrence equation, it is possible to determine a basis of its
solution space in terms of asymptotic series, see [WZ], [K11], [K23] for how
this is done. 

For example, the algebraic power series
$$
  F(X) = \frac1{\sqrt{1-2X}}
$$
satisfies the differential equation
$$
  (1-2X)\,F'(X) - F(X) = 0\quad,
$$
and its coefficient sequence, $f_n$, satisfies the recurrence
$$
  (n+1)\,f_{n+1} - (2n+1)\,f_n = 0\quad.
$$
This recurrence has the asymptotic series solution
$$
  2^n \, n^{-1/2} \, (1 - \frac18n^{-1} + \frac1{128}n^{-2} + \frac5{1024}\,n^{-3} + \cdots)\quad,
$$
and, since the recurrence has order one, the space of all solutions consists of all
constant multiples of this series.

In general, an algebraic power series leads to a recurrence of higher order, and in that case,
there will be several linearly independent asymptotic series solutions. Fortunately, in many
cases (at least in all cases considered in this paper), one of the series solutions dominates
all the others, and the asymptotics of the coefficient series is simply a constant multiple
of the series solution with the fastest growth rate.

So for a given algebraic power series $F(X)$ whose coefficient series grows like
$$
  C\,\mu^n\,n^\theta(1 + \frac{C_1}{n} + \frac{C_2}{n^2} + \cdots)\quad,
$$
we can find $\mu$, $\theta$, and $C_1,C_2,\dots$ by first translating the algebraic equation
into a differential equation, then translating the differential equation into a recurrence,
and then computing the dominating asymptotic series solution of this recurrence.

It remains to determine the multiplicative constant~$C$.

For sequences defined by arbitrary linear recurrences with polynomial
coefficients, it is not known in general how to identify this constant exactly,
but since the coefficient sequences of algebraic power series are more special, 
this can be done. Following the advice of [FS], we should view the power series as an
analytic function rather than merely as a formal object. The asymptotics of the
series coefficients at the origin are then determined by the local behaviour of
the function near the singularity that is closest to the origin.

More precisely, if $\zeta$ is the singularity closest to the origin, and
near this singularity (approached from the right direction) we have
$$
  F(X) \sim A\, (1 - \frac{X}{\zeta})^{k/r}
$$
for some rational number $k/r$ (not a positive integer), then the asymptotics
of the series coefficients of $F(X)$ is given by
$$
  f_n \sim \frac{A}{\Gamma(-k/r)}\, \zeta^{-n} \, n^{-1-k/r}\qquad(n\to\infty)\quad.
$$
Once more, we refer to [FS] for the full story.

In the example discussed earlier, we have $\zeta=1/2$ and
$$
  F(X) = \frac1{\sqrt{1-2X}} = 1\,(1-\frac{X}{1/2})^{-1/2}\quad,
$$
which confirms the terms $2^n$ and $n^{-1/2}$ that we found before, and in
addition reveals $C=1/\Gamma(1/2)=1/\sqrt{\pi}$. 

Every singularity of an algebraic function is either a pole or a branch point,
and both can be easily computed from the minimal polynomial of the function. It
is also easy to determine which of them is closest to the origin, at least for
algebraic functions of degree two. For higher degrees, things get more delicate
because we must be sure not to be working on the wrong branch of the
function. We do not need to worry about this issue here because we are only
concerned with algebraic functions of degree two.

Once the nearest singularity to the origin has been identified, the required
data can be obtained by computing the first term of the series expansion of
the algebraic function at that point.
All these computations can be performed with exact algebraic numbers, thus
establishing a rigorous proof of Theorem~3. The Mathematica code employed
for these computations is here:

{\tt http://sites.math.rutgers.edu/\~{}zeilberg/tokhniot/AlgAsympt.m} \quad .

{\bf References}

[B] AJ Bu, {\it A Maple package for computing generating functions enumerating RNA secondary structures}.
The Maple package is here: \hfill \break
{\tt https://sites.math.rutgers.edu/\~{}zeilberg/tokhniot/RNA.txt} \quad . \hfill \break
Here is the input file: {\tt https://sites.math.rutgers.edu/\~{}zeilberg/tokhniot/inRNA1.txt} \quad , \hfill \break
that produces the following output file:\hfill\break
{\tt https://sites.math.rutgers.edu/\~{}zeilberg/tokhniot/oRNA1.txt} \quad .

[CRU] Sang Kwan Choi, Chaiho Rim, and Hwajin Um, {\it Narayana number, Chebyshev polynomial, and Motzkin path on RNA abstract shapes},
in: {\it 2017 Matrix Annals}, edited by David R. Wood et. al, Springer, 2019. \hfill \break
{\tt https://sites.math.rutgers.edu/\~{}zeilberg/akherim/choi2019.pdf} \quad .

[FS] Phillipe Flajolet and Robert Sedgewick, {``Analytic Combinatorics''}, Cambridge University Press, 2009. Freely available on-line: \hfill\break
{\tt https://ac.cs.princeton.edu/home/AC.pdf} \quad .

[GD] Shile Gao, and Kequan Ding, {\it A graphical criterion of planarity for RNA secondary structures with pseudoknots in Rivas-Eddy class},
Theoretical Computer Science {\bf 395} (2008), 47-56.  \hfill \break
{\tt https://www.sciencedirect.com/science/article/pii/S0304397507006858} \quad .

[K11] Manuel Kauers, {\it A {M}athematica Package for Computing Asymptotic Expansions of Solutions of P-Finite Recurrence Equations},
Technical Report 11-04, RISC-Linz. (2011).

[K23] Manuel Kauers, {\it D-Finite Functions}, Springer, 2023.

[KP] Manuel Kauers and Peter Paule, {\it The Concrete Tetrahedron}, Springer, 2011.

[Sl] Neil A. J. Sloane, {\it The On-Line Encyclopedia of Integer Sequences® (OEIS)}, {\tt https://oeis.org/}.

[W] M.S. Waterman, {\it Secondary structure of single-stranded nucleic acids}, Adv. Math. Supp.  Stud. {\bf 1} (1978), 167--212. \hfill \break
{\tt https://dornsife.usc.edu/msw/wp-content/uploads/sites/236/2023/09/msw-026.pdf} \quad .

[WZ] Jet Wimp and Doron Zeilberger, {\it Resurrecting the Asymptotics of Linear Recurrences},
J. Math. Analysis and Appl. {\bf 111} (1985), 162--176. \hfill\break
{\tt https://sites.math.rutgers.edu/\~{}zeilberg/mamarimY/WimpZeilberger.pdf}.

\bigskip
\hrule
\bigskip
AJ Bu,  Department of Mathematics, Rutgers University (New Brunswick), Hill Center-Busch Campus, 110 Frelinghuysen
Rd., Piscataway, NJ 08854-8019, USA. \hfill\break
Email: {\tt   ab1854  at math dot rutgers dot edu}
\bigskip
Manuel Kauers, Institute for Algebra, Johannes Kepler Universit\"at, Altenbergerstra\ss e 69 A-4040 Linz, Austria
Email: {\tt   manuel dot kauers at jku dot at}
\bigskip
Doron Zeilberger,  Department of Mathematics, Rutgers University (New Brunswick), Hill Center-Busch Campus, 110 Frelinghuysen
Rd., Piscataway, NJ 08854-8019, USA. \hfill\break
Email: {\tt   ab1854  at math dot rutgers dot edu}

Written: {\bf  Feb. 20, 2026}.

\end